\documentclass[a4paper, leqno, 12pt]{article}
\usepackage{amsmath,amsthm}
\usepackage{amssymb,latexsym}
\usepackage{enumerate}
\usepackage{color}
\usepackage[T1]{fontenc}
\usepackage{array}
\usepackage[all]{xypic}

\usepackage[left=3cm,right=3cm,top=3.0cm,bottom=3.0cm]{geometry}

\newtheorem{thm}{Theorem}[section]
\newtheorem{cor}[thm]{Corollary}
\newtheorem{lem}[thm]{Lemma}
\newtheorem{prop}[thm]{Proposition}

\theoremstyle{definition}

\numberwithin{equation}{section}

 1

\DeclareMathOperator{\precc}{\preccurlyeq}

\DeclareMathOperator{\LM}{\textnormal{LM}}\DeclareMathOperator{\LT}{\textnormal{LT}}

\DeclareMathOperator{\N}{\mathbb{N}}\DeclareMathOperator{\T}{\mathbb{T}}

\DeclareMathOperator{\ra}{\rightarrow}

\DeclareMathOperator{\lra}{\leftrightarrow}
 \DeclareMathOperator{\rag}{\stackrel{\textit{G}}{\rightarrow}}
\DeclareMathOperator{\ran}{\rangle}\DeclareMathOperator{\lan}{\langle}
\DeclareMathOperator{\para}{\parallel}\DeclareMathOperator{\npara}{\nparallel}

\def\NN{{\mathbb{N}}}
\def\TT{{\mathbb{T}}}
\def\FF{{\mathbb{F}}}

\def\CB{{\cal B}}

\def\CF{{\cal F}}

\def\epv {{$\mbox{}$\hfill ${\Box}$\vspace*{1.5ex} }}

\def\lm{\textnormal{lm}} \def\lt{\textnormal{lt}}

\def\ra{\rightarrow}

\def\und#1{\underline{#1}}

\begin{document}


\baselineskip=17pt


\title{On ascending chains of ideals in the polynomial ring}
\author{Grzegorz Pastuszak (Toru\'n)}

\date{}

\maketitle

\renewcommand{\thefootnote}{}

\renewcommand{\thefootnote}{\arabic{footnote}}
\setcounter{footnote}{0}


\begin{abstract}
Assume that $K$ is a field and $I_{1}\subsetneq ...\subsetneq I_{t}$ is an ascending chain (of length $t$) of ideals in the polynomial ring $K[x_{1},,...,x_{m}]$, for some $m\geq 1$. Suppose that $I_{j}$ is generated by polynomials of degrees less or equal to some natural number $f(j)\geq 1$, for any $j=1,...,t$. In the paper we construct, in an elementary way, a natural number $\CB(m,f)$ (depending on $m$ and the function $f$) such that $t\leq\CB(m,f)$. We also discuss some possible applications of this result.
\end{abstract}

\section{Introduction}

Assume that $K$ is a field and $K[x_{1},,...,x_{m}]$ is the polynomial ring over $K$ in $m\geq 1$ variables. Denote by $\NN_{1}$ the set of all natural numbers greater or equal to $1$ and let $f:\NN_{1}\ra\NN_{1}$ be an arbitrary function. Assume that $$I_{1}\subsetneq ...\subsetneq I_{t}\subseteq K[x_{1},,...,x_{m}]$$ is an ascending chain (of length $t$) of ideals such that $I_{j}$ is generated by polynomials of degrees less or equal to $f(j)$, for any $j=1,...,t$. 

In \cite{Se1} A. Seidenberg shows that there exists a natural number $g_{m}(f)$, for an increasing $f$, such that $t\leq g_{m}(f)$. He proposes rather complicated, but an explicit formula for $g_{m}(f)$ in terms of $m$ and $f$. In \cite{Mo} G. Moreno Soc{\'i}as finds a better bound for the number $t$ and expresses it, in terms of $m$ and $f$, in a quite optimal way. He also shows, among other things, that the number $g_{m}(f)$ is primitive recursive in $f$, for any $m\geq 1$. Another approach to the problem is given in \cite{AsPo} where the authors obtain more general facts in somewhat extended context. For example, Proposition 3.22 from \cite{AsPo} implies some of the main results of \cite{Se1} and \cite{Mo}. Note that both \cite{Mo} and \cite{AsPo} widely use the \textit{Hilbert-Samuel polynomials} and related concepts, see for example \cite[Chapter 4]{BH} and \cite[Section 19.5]{E}.

This paper is devoted to construct the number $g_{m}(f)$, denoted here by $\CB(m,f)$, in an elementary way. We apply only some basic facts from the theory of Gr\"obner bases. 

The paper is organized as follows. In Section 2 we fix the notation and recall some information about Gr\"obner bases, e.g. the renowned algorithm for constructing a Gr\"obner basis of a given ideal, due to B. Buchberger. 

Section 3 is the core of the paper. In Theorem 3.5 (concluding all the preceding results) we define a function $\CB$ with the \textit{bounding property} which sets a bound on the length of antichains in $\NN^{m}$, see Sections 2 and 3 for all the definitions. Our arguments are combinatorial and rather elementary. Theorem 3.5 is further applied in the next section.

In Section 4 we present the main results of the paper. We show how to reduce the general problem studied in the paper to the situation considered in Section 3. The main result on ascending chains of ideals in $K[x_{1},...,x_{m}]$ is given in Theorem 4.2. Furthermore, we derive some interesting consequences of Theorem 4.2 in Corollaries 4.4 and 4.5. 

In the last section of the paper we describe our motivation to study bounds of ascending chains of ideals in the polynomial ring. As we write in detail in Section 5, the motivation comes from the first order logic and elimination of quantifiers. Namely, in the subsequent paper \cite{P} we apply Corollary 4.5 to give a constructive proof of Tarski's theorem on quantifier elimination in the theory of algebraically closed fields. In a sense, the present paper rediscovers some of the main results of \cite{Mo} and \cite{AsPo} in order to prove Tarski's theorem in a constructive way. 

The results presented in the paper are part of the author's master's thesis, supervised by Stanis{\l}aw Kasjan in 2007. The author is grateful to the supervisor for all discussions and support during the work on the thesis.

\section{Gr\"obner bases and Buchberger's algorithm}


We denote by $\NN$ the set of all natural numbers and by $\NN_{1}$ the set $\NN\setminus\{0\}$. Assume that $K$ is a field and $m\in\NN_{1}$. Then $K[x_{1},...,x_{m}]$ is the polynomial ring over $K$ in $m$ variables $x_{1},...,x_{m}$. The set of all monomials in $K[x_{1},...,x_{m}]$ is denoted by $\TT_{m}$. If $\alpha=(a_{1},...,a_{m})\in\NN^{m}$, then the monomial $x_{1}^{a_{1}}...x_{m}^{a_{m}}\in\TT_{m}$ is denoted by $\und{x}^{\alpha}$. The \textit{degree} of $\und{x}^{\alpha}=x_{1}^{a_{1}}...x_{m}^{a_{m}}$ is the sum $a_{1}+...+a_{m}$. A polynomial $f\in K[x_{1},...,x_{m}]$ is denoted by $\sum_{\alpha}a_{\alpha}\und{x}^{\alpha}$ where $a_{\alpha}\in K$ and $a_{\alpha}=0$ for almost all $\alpha\in\N^{m}$. If $f=\sum_{\alpha}a_{\alpha}\und{x}^{\alpha}$, then the set $\{\und{x}^{\alpha};a_{\alpha}\neq 0\}$ is the \textit{support} of $f$. The \textit{degree} of $f$, denoted by $\deg(f)$, is the maximum of degrees of monomials from the support of $f$.

Assume that $m\in\NN_{1}$. We view the set $\NN^{m}$ as a monoid with respect to the pointwise addition, denoted by $+$. We denote by $\und{0}$ the neutral element $(0,...,0)\in\NN^{m}$ of $+$. If $\alpha,\beta\in\NN^{m}$ and $\alpha+\gamma=\beta$ for some $\gamma\in\NN^{m}$, then we write $\alpha\para\beta$. Note that $\para$ defines an order on $\NN^{m}$ and $\NN^{m}$ is an ordered monoid with respect to $+$ and $\para$. Obviously, $\alpha\para\beta$ if and only if $\und{x}^{\alpha}$ divides $\und{x}^{\beta}$. If $\alpha\in\NN^{m}$ and $\alpha=(a_{1},...,a_{m})$, then we set $|\alpha|=a_{1}+...+a_{m}$ and hence $\deg(\und{x}^{\alpha})=|\alpha|$. Recall that a binary relation $\precc$ on $\NN^{m}$ is an \textit{admissible relation} (or an \textit{admissible ordering}) if and only if the following three conditions are satisfied: $\precc$ is a linear ordering, $\und{0}\precc\alpha$ for any $\alpha\in\NN^{m}$ and $\alpha\precc\beta$ yields $\alpha+\gamma\precc\beta+\gamma$ for any $\alpha,\beta,\gamma\in\N^{m}$. Note that $\alpha\para\beta$ implies $\alpha\precc\beta$ and any admissible relation is a well-order, see Chapter 1 of \cite{AL}. We call an admissible relation $\precc$ on $\NN^{m}$ \textit{graded} if and only if $\alpha\precc\beta$ implies $|\alpha|\leq|\beta|$ for any $\alpha,\beta\in\NN^{m}$. A basic example of an admissible relation is the \textit{lexicographical order}. Its graded version is called the \textit{degree lexicographical order}. We send to \cite{AL} for definitions of these orders, as well as for other examples.

It is easy to see that an admissible relation on $\NN^{m}$ induces a relation on the set $\TT_{m}$ of all monomials in $K[x_{1},...,x_{m}]$ via the natural identification $(a_{1},...,a_{m})\lra x_{1}^{a_{1}}...x_{m}^{a_{m}}$. We call such a relation a \textit{monomial ordering}.

Assume that $\precc$ is an admissible relation on $\NN^{m}$. If $f=\sum_{\alpha}a_{\alpha}\und{x}^{\alpha}$ and $\eta$ is the greatest element of the set $\{\alpha\in\NN^{m};a_{\alpha}\neq 0\}$ with respect to $\precc$, then $\und{x}^{\eta}$ is the \textit{leading monomial} of $f$ (denoted by $\lm(f)$) and $a_{\eta}\und{x}^{\eta}$ is the \textit{leading term} of $f$ (denoted by $\lt(f)$). If $I$ is an ideal in $K[x_{1},...,x_{m}]$, then we set $\LM(I)=\{\lm(f);f\in I\}$ and $\LT(I)=\{\lt(f);f\in I\}$. 

Assume that $f,f_{1},...,f_{s}\in K[x_{1},...,x_{m}]$ and set $F=\{f_{1},...,f_{s}\}$. Then there are $a_{1},...,a_{s},r\in K[x_{1},...,x_{m}]$ such that $f=a_{1}f_{1}+...+a_{s}f_{s}+r$, $\lm(f)$ is the greatest element of the set $\{\lm(a_{1}f_{1}),...,\lm(a_{s}f_{s}),\lm(r)\}$ and $r$ is \textit{reduced} modulo $F$, that is, $\lm(f_{i})$ does not divide any element of support of $r$, for any $i=1,...,s$. In this case we say that $r$ is a \textit{reduction} of $f$ modulo $F$ and we write $f\stackrel{F}{\ra}r$ or $r=f_{F}$. A reduction $r$ of $f$ modulo $F$ is the result of the \textit{Multivariable Division Algorithm}, see for example \cite[I.5]{AL}.

Assume that $I$ is an ideal in $K[x_{1},...,x_{m}]$ and $\precc$ is an admissible relation on $\NN^{m}$. A set $G=\{g_{1},...,g_{t}\}\subseteq I$ is a \textit{Gr\"obner basis} of $I$ (with respect to $\precc$) if and only if, for any $f\in I$, there is $i=1,...,t$ such that $\lm(g_{i})$ divides $\lm(f)$.

For the rest of the section $\precc$ denotes a fixed admissible relation on $\NN^{m}$. The following theorem is a basic result in the theory of Gr\"obner bases.

\begin{thm} Assume that $I$ is a non-zero ideal in $K[x_{1},...,x_{m}]$ and $G=\{g_{1},...,g_{t}\}$, $G\subseteq I$, is a set of non-zero polynomials. The following conditions are equivalent.
\begin{enumerate}[\rm(1)]
	\item The set $G$ is a Gr\"obner basis of $I$.
	\item $f\in I$ if and only if $f\rag 0$.
	\item $f\in I$ if and only if there are polynomials $h_{1},...,h_{t}$ such that $f=\sum_{i=1}^{t}h_{i}g_{i}$ and $\lm(f)=\max\{\lm(h_{1}g_{1}),...,\lm(h_{1}g_{1})\}$.
	\item $\lan LM(G)\ran=\lan LM(I)\ran$.
\end{enumerate}
\end{thm}

{\bf Proof.} See the proof of \cite[Theorem 1.6.2]{AL}. \epv

The above theorem yields that if $G$ is a Gr\"obner basis of $I$, then $\lan G \ran=I$. Hence we say that a finite set of polynomials $G$ is a \textit{Gr\"obner basis} if and only if $G$ is a Gr\"obner basis of $\lan G \ran$. Theorem 2.1 also implies that any non-zero ideal in $K[x_{1},...,x_{m}]$ has a Gr\"obner basis.

The definition of Gr\"obner basis was introduced by B. Buchberger in \cite{B1}. Now we present a fundamental method for constructing a Gr\"obner basis of a given ideal, known as the \textit{Buchberger's algorithm}, which is also given in \cite{B1}. We start with the following crucial notion of $S$-polynomial.

Assume that $f,g\in K[x_{1},...,x_{m}]$, $f,g\neq 0$ and $\und{x}^{\alpha}$ is the greatest common multiple of $\lm(f)$ and $\lm(g)$. Then the polynomial $$S(f,g)=\frac{\und{x}^{\alpha}}{\lt(f)}f-\frac{\und{x}^{\alpha}}{\lt(g)}g$$ is the \textit{$S$-polynomial} of $f$ and $g$. If $B=\{b_{1},...,b_{s}\}$ is a finite set of polynomials in $K[x_{1},...,x_{m}]$, then we define $S_{B}$ to be the set of all non-trivial reductions of $S$-polynomials of $b_{i}$ and $b_{j}$ modulo $B$, that is, $$S_{B}=\{S(b_{i},b_{j})_{B};b_{i},b_{j}\in B\}\setminus\{0\}.$$ 

The following fact from \cite{B1} (see also \cite{B2}) sets the ground for the succeeding Buchberger's algorithm.

\begin{thm} Assume that $G=\{g_{1},...,g_{t}\}$ is a set of non-zero polynomials in $K[x_{1},...,x_{m}]$. Then $G$ is a Gr\"obner basis if and only if $S(g_{i},g_{j})\rag 0$ for any $i,j$.
\end{thm}

{\bf Proof.} See the proof of \cite[Theorem 1.7.4]{AL}. \epv

\textbf{Algorithm (B. Buchberger).} Input: a set $F=\{f_{1},...,f_{s}\}\subseteq K[x_{1},...,x_{n}]$ of non-zero polynomials. Output: a set $G=\{g_{1},...,g_{t}\}\subseteq K[x_{1},...,x_{n}]$ such that $F\subseteq G$ and $G$ is a Gr\"obner basis of $\lan F\ran$.

\begin{enumerate}[\rm (1)]
	\item Set $B_{0}:=F$ and $i:=0$.
  \item Put $B_{i+1}:=B_{i}\cup S_{B_{i}}$. If $B_{i+1}\neq B_{i}$, then put $i:=i+1$ and return to (2). Otherwise put $G:=B_{i}$ and finish.
\end{enumerate}

Theorem 2.2 yields that the Buchberger's algorithm is correct. Note that this algorithm halts, because $\lan\LT(B_{i})\ran\subsetneq\lan\LT(B_{i+1})\ran$ for any $i\geq 0$ and, in a noetherian ring, any ascending chain of ideals is finite.

\section{Antichains in $\NN^{m}$}

A sequence $\alpha_{1},...,\alpha_{t}\in\NN^{m}$ is an \textit{antichain} if and only if $\alpha_{i}\npara\alpha_{j}$ for any $i<j$. Denote by $\FF$ the set of all non-decreasing functions $\N_{1}\ra\N_{1}$ and let $f\in\FF$. We say that an antichain $\alpha_{1},...,\alpha_{t}\in\NN^{m}$ is \textit{$f$-bounded} if and only if $|\alpha_{i}|\leq f(i)$ for any $i=1,...,t$. 

In this section we give a bound on the length of $f$-bounded antichains in $\NN^{m}$ depending on $m\in\NN_{1}$ and $f\in\FF$. Let us start with some notation and terminology.

We write $f\leq f'$ if and only if $f(n)\leq f'(n)$ for any $n\in\NN_{1}$ and $f,f'\in\FF$. Assume that $m\geq 1$ is a natural number. We say that a function $\CB_{m}:\FF\ra\NN$ has the \textit{bounding property for $m$} if and only if the following conditions are satisfied:
\begin{enumerate}[\rm(1)]
	\item $t\leq\CB_{m}(f)$ for any $f\in\FF$ and $f$-bounded antichain $\alpha_{1},...,\alpha_{t}\in\NN^{m}$ of length $t$,
	\item $\CB_{m}(f)\leq\CB_{m}(f')$ for any $f,f'\in\FF$ such that $f\leq f'$.
\end{enumerate}

We say that a function $\CB:\N_{1}\times\FF\ra\NN$ has the \textit{bounding property} if and only if, for any $m\in\N_{1}$, the function $\CB_{m}:\FF\ra\NN$ defined by $\CB_{m}(f)=\CB(m,f)$, for any $f\in\FF$, has the bounding property for $m$.  

This section is devoted to construct a function with the bounding property in the above sense. As an equivalent, we construct a sequence $(\CB_{m})_{m\in\NN_{1}}$ of functions such that $\CB_{m}$ has the bounding property for $m$. Our construction is inductive with respect to the number $m$.

The existence of a function with the bounding property is rather straightforward consequence of the Compactness Theorem of first order logic, see \cite{F} and \cite[Proposition 3.25]{AsPo} for more details. However, this approach does not provide the explicit form of a function with the bounding property.

In the following proposition we construct a function $\CB_{1}:\FF\ra\NN$ with the bounding property for $m=1$. This is the first step of our induction. 

\begin{prop} The function $\CB_{1}:\FF\ra\NN$ such that $\CB_{1}(f)=f(1)+1$, for any $f\in\FF$, has the bounding property for $m=1$.
\end{prop}

{\bf{Proof.}} Assume that $f\in\FF$ and $\alpha_{1},...,\alpha_{t}\in\NN$ is an $f$-bounded antichain. Then $f(1)\geq\alpha_{1}>\alpha_{2}>...>\alpha_{t}$ and so $t\leq f(1)+1=\CB_{1}(f)$. Moreover, if $f,g\in\N_{1}$ and $f\leq g$, then $\CB_{1}(f)=f(1)+1\leq g(1)+1=\CB_{1}(g)$. This yields $\CB_{1}:\FF\ra\NN$ has the bounding property for $m=1$. \epv

Before the second step of the induction, we introduce the following terminology which generalizes, in some sense, the one given before.



Assume that $m\geq 1$, $\alpha_{1}=(a_{11},a_{12},...,a_{1m}),...,\alpha_{t}=(a_{t1},a_{t2},...,a_{tm})\in\NN^{m}$ is an antichain, $\beta=(b_{1},...,b_{k})\in\NN^{k}$, for some $k\in\{1,...,m\}$ (we treat $\beta$ as the sequence $b_{1},...,b_{k}$), and $f\in\FF$. We say that the chain $\alpha_{1},...,\alpha_{t}$ is \textit{$(f,\beta)$-bounded} (or \textit{$(f,b_{1},...,b_{k})$-bounded}) if and only if it is $f$-bounded and $$a_{11},a_{21},...,a_{t1}\leq b_{1},$$ $$a_{12},a_{22},...,a_{t2}\leq b_{2},$$ $$\vdots$$ $$a_{1k},a_{2k},...,a_{tk}\leq b_{k}.$$

We say that a function $\CB_{m}^{k}:\FF\times\NN^{k}\ra\NN$ has the \textit{$k$-bounding property for $m$} if and only if the following conditions are satisfied:
\begin{enumerate}[\rm(1)]
	\item $t\leq\CB_{m}^{k}(f,\beta)$ for any $f\in\FF$, $\beta\in\NN^{k}$ and $(f,\beta)$-bounded antichain $\alpha_{1},...,\alpha_{t}\in\NN^{m}$ of length $t$,
	\item $\CB_{m}^{k}(f,\beta)\leq\CB_{m}^{k}(f',\beta')$ for any $f,f'\in\FF$ and $\beta,\beta'\in\NN^{k}$ such that $f\leq f'$ and $\beta\para\beta'$.
\end{enumerate}

Recall that if $\beta=(b_{1},...,b_{k})$ and $\beta'=(b'_{1},...,b'_{k})$, then the condition $\beta\para\beta'$ means $b_{i}\leq b'_{i}$ for any $i=1,...,k$.

We agree that a function $\CB_{m}:\FF\ra\NN$ with the bounding property for $m$ has the $0$-bounding property for $m$ (and vice versa).

Assume that the function $\CB_{m-1}:\FF\ra\NN$, $m\geq 2$, has the bounding property for $m-1$. Our aim is to construct a function $\CB_{m}:\FF\ra\NN$ with the bounding property for $m$. In order to do this, we construct functions $\CB_{m}^{k}:\FF\times\NN^{k}\ra\NN$ having $k$-bounding properties for $m$ by the backward induction with respect to $k$. To be more precise, we first construct the function $\CB_{m}^{m}:\FF\times\NN^{m}\ra\NN$ having the $m$-bounding property for $m$ (this construction is general and does not depend on $\CB_{m-1}:\FF\ra\NN$, see Proposition 3.2). Then we show how to obtain $\CB_{m}^{k}:\FF\times\NN^{k}\ra\NN$ from $\CB_{m}^{k+1}:\FF\times\NN^{k+1}\ra\NN$. This process provides a function with the $0$-bounding property for $m$, that is, a function with the bounding property for $m$.

The first step of the backward induction is given in the following fact. 

\begin{prop} Assume that $m\geq 1$. The function $\CB_{m}^{m}:\FF\times\NN^{m}\ra\NN$ such that $\CB_{m}^{m}(f,b_{1},...,b_{m})=(b_{1}+1)\cdot...\cdot(b_{m}+1)$ has the $m$-bounding property for $m$.
\end{prop}

{\bf{Proof.}} Assume that $f\in\FF$ and $b_{1},...,b_{m}\in\NN$. The set of all $m$-tuples $(a_{1},...,a_{m})$ of natural numbers such that $a_{i}\leq b_{i}$, for $i=1,...,m$, has $(b_{1}+1)\cdot...\cdot(b_{m}+1)$ elements. This shows that if $\alpha_{1},...,\alpha_{t}\in\NN$ is an $(f,b_{1},...,b_{m})$-bounded antichain, then $$t\leq (b_{1}+1)\cdot...\cdot(b_{m}+1)=\CB_{m}^{m}(f,b_{1},...,b_{m}).$$ Moreover, if $b_{i}\leq b'_{i}$ for $i=1,...,m$, then $\CB_{m}^{m}(f,b_{1},...,b_{m})\leq\CB_{m}^{m}(g,b'_{1},...,b'_{m})$ for any $f,g\in\FF$. Hence $\CB_{m}^{m}:\FF\times\NN^{m}\ra\NN$ has the $m$-bounding property for $m$. \epv

Now we introduce some notation. If $\alpha=(a_{1},...,a_{m})\in\NN^{m}$ and $s\in\{1,...,m\}$, then we set $\widehat{\alpha}^{s}=(a_{1},...,a_{s-1},a_{s+1},...,a_{m})\in\NN^{m-1}$.

If $f\in\FF$ and $s\in\NN$, then ${}^{s}f:\NN_{1}\ra\NN_{1}$ is a function such that ${}^{s}f(n)=f(s+n)$ for any $n\in\NN_{1}$. Observe that ${}^{s}f\in\FF$.

Assume that $m\geq 2$, $k\in\{0,...,m-1\}$ and the function $\CB_{m}^{k+1}:\FF\times\NN^{k+1}\ra\NN$ has the $(k+1)$-bounding property for $m$. Suppose $f\in\FF$, $\beta\in\NN^{k}$ and define recursively a function $g:\NN_{1}\ra\NN_{1}$ in the following way:
\begin{enumerate}[\rm(1)]
	\item $g(1)=1$,
	\item $g(n+1)=1+g(n)+\CB^{k+1}_{m}({}^{g(n)}f,\beta,f(g(n)))$ for any $n\geq 1$.
\end{enumerate}
Obviously $g\in\FF$ and hence we get a function $\CF_{m}^{k}:\FF\times\NN^{k}\ra\FF$ such that $(f,\beta)\mapsto g$. We use this function in the following lemma which is the key ingredient of the second step of the backward induction.

\begin{lem} Assume that $m\geq 2$, $k\in\{0,...,m-1\}$ and $\CB_{m}^{k+1}:\FF\times\NN^{k+1}\ra\NN$ has the $(k+1)$-bounding property for $m$. Assume that $f\in\FF$, $\beta\in\NN^{k}$ and $\alpha_{1},...,\alpha_{t}$ is an $(f,\beta)$-bounded antichain in $\NN^{m}$ of length $t$.

\begin{enumerate}[\rm (1)]
	\item Assume that $\alpha_{s_{1}},...,\alpha_{s_{r}}$ is a subsequence of $\alpha_{1},...,\alpha_{t}$ such that the sequence $\widehat{\alpha}_{s_{1}}^{k+1},...,\widehat{\alpha}_{s_{r}}^{k+1}\in\NN^{m-1}$ is an antichain. Set $\mu=s_{r}+\CB_{m}^{k+1}({}^{s_{r}}f,\beta,f(s_{r}))$. If we have $\mu+1\leq t$, then there is a natural number $c\in\{s_{r}+1,...,\mu+1\}$ such that the sequence $\widehat{\alpha}_{s_{1}}^{k+1},...,\widehat{\alpha}_{s_{r}}^{k+1},\widehat{\alpha}_{c}^{k+1}$ is an antichain in $\NN^{m-1}$.
	\item Set $g=\CF_{m}^{k}(f,\beta)$, fix a natural number $n\geq 1$ and suppose that $g(n)\leq t$. Then there is a subsequence $\alpha_{p_{1}},...,\alpha_{p_{n}}$ of length $n$ of the sequence $\alpha_{1},...,\alpha_{t}$ such that the sequence $\widehat{\alpha}_{p_{1}}^{k+1},...,\widehat{\alpha}_{p_{n}}^{k+1}$ is an $(f\circ g)$-bounded antichain in $\NN^{m-1}$.
\end{enumerate}
\end{lem}

{\bf{Proof.}} (1) Set $\alpha_{1}=(a_{11},a_{12},...,a_{1m}),...,\alpha_{t}=(a_{t1},a_{t2},...,a_{tm})$ and assume $\mu+1\leq t$, $d\geq s_{r}+1$. Suppose that, for any $n\in\{s_{r}+1,...,d\}$, the sequence $\widehat{\alpha}_{s_{1}}^{k+1},...,\widehat{\alpha}_{s_{r}}^{k+1},\widehat{\alpha}_{n}^{k+1}$ is not an antichain in $\NN^{m-1}$. We show that $d\leq\mu$. Indeed, for a fixed $n$ we have $\widehat{\alpha}_{s_{i}}^{k+1}||\widehat{\alpha}_{n}^{k+1}$ for some $i$, because $\widehat{\alpha}_{s_{1}}^{k+1},...,\widehat{\alpha}_{s_{r}}^{k+1}$ is an antichain in $\NN^{m-1}$. Note that 
$\alpha_{s_{1}},...,\alpha_{s_{r}},\alpha_{n}$ is an antichain in $\NN^{m}$, so $\alpha_{s_{i}}\npara\alpha_{n}$. Hence we get $a_{s_{i}(k+1)}>a_{n(k+1)}$ and $f(s_{r})\geq f(s_{i})\geq a_{s_{i}(k+1)}>a_{n(k+1)}$.

Consequently, $a_{n(k+1)}\leq f(s_{r})$ for any $n\in\{s_{r}+1,...,d\}$ and thus the sequence $\alpha_{s_{r}+1},\alpha_{s_{r}+2},...,\alpha_{d}\in\NN^{m}$ is $({}^{s_{r}}f,\beta,f(s_{r}))$-bounded. This implies that $d-s_{r}\leq\CB_{m}^{k+1}({}^{s_{r}}f,\beta,f(s_{r}))$, so $d\leq\mu$ and $(1)$ follows. 

(2) We use induction with respect to $n$. Assume that $n=1$ and set $p_{1}=1$. Then $\widehat{\alpha}_{1}^{k+1}$ is an $(f\circ g)$-bounded antichain in $\NN^{m-1}$, because $|\widehat{\alpha}_{1}^{k+1}|\leq f(1)=f(g(1))$.

Assume that the thesis holds for some $n\geq 1$. Moreover, assume a technical condition $p_{1}\leq g(1),...,p_{n}\leq g(n)$. We show that the thesis holds for $n+1$ and $p_{1}\leq g(1),...,p_{n+1}\leq g(n+1)$. Indeed, if $g(n+1)\leq t$, then $g(n)\leq t$ and hence there is an antichain in $\NN^{m-1}$ of the form $\widehat{\alpha}_{p_{1}}^{k+1},...,\widehat{\alpha}_{p_{n}}^{k+1}.$ Observe that $$g(n+1)=1+g(n)+\CB^{k+1}_{m}({}^{g(n)}f,\beta,f(g(n)))\leq t$$ and thus, applying $(1)$ for $s_{r}=g(n)$, we get $c\in\{g(n)+1,...,g(n+1)\}$ such that the sequence $\widehat{\alpha}_{p_{1}}^{k+1},...,\widehat{\alpha}_{p_{n}}^{k+1},\widehat{\alpha}_{c}^{k+1}$ is an antichain in $\NN^{m-1}$. Set $p(n+1)=c$. Since $p_{n+1}\leq g(n+1)$, we get $|\widehat{\alpha}_{p_{n+1}}^{k+1}|\leq f(p_{n+1})\leq f(g(n+1))$ and thus this antichain is $(f\circ g)$-bounded. This finishes the proof. \epv

Given the above lemma we are able to prove the second step of the backward induction with respect to $k$ and hence the second step of the main induction (with respect to $m$).

\begin{cor} Assume $m\geq 2$, $k\in\{0,...,m-1\}$, $\CB_{m}^{k+1}:\FF\times\NN^{k+1}\ra\NN$ has the $(k+1)$-bounding property for $m$ and $\CB_{m-1}:\FF\ra\NN$ has the bounding property for $m-1$. 
\begin{enumerate}[\rm (1)]
	\item The function $\CB_{m}^{k}:\FF\times\NN^{k}\ra\NN$ such that $$\CB_{m}^{k}(f,\beta)=g(\CB_{m-1}(f\circ g)+1),$$ for any $f\in\FF$, $\beta\in\NN^{k}$ and $g=\CF_{m}^{k}(f,\beta)$, has the $k$-bounding property for $m$.
	\item The function $\CB_{m}^{0}:\FF\ra\NN$ has the bounding property for $m$.
\end{enumerate}
\end{cor}

{\bf Proof.} (1) Assume that the antichain $\alpha_{1},...,\alpha_{t}\in\NN^{m}$ of length $t$ is $(f,\beta)$-bounded. If $g(\CB_{m-1}(f\circ g)+1)\leq t$, then Lemma 3.3 (2) implies there is an $(f\circ g)$-bounded antichain in $\NN^{m-1}$ of length $\CB_{m-1}(f\circ g)+1$, a contradiction. Hence $t<g(\CB_{m-1}(f\circ g)+1)$ and it is enough to prove that $\CB_{m}^{k}(f,\beta)\leq\CB_{m}^{k}(f',\beta')$ for any $f,f'\in\FF$ and $\beta,\beta'\in\NN^{k}$ such that $f\leq f'$ and $\beta||\beta'$. Set $g=\CF_{m}^{k}(f,\beta)$ and $g'=\CF_{m}^{k}(f',\beta')$. It follows easily from the construction of $g,g'$ that $g\leq g'$. Thus $f\circ g\leq f'\circ g'$, $\CB_{m-1}(f\circ g)\leq\CB_{m-1}(f'\circ g')$ and finally $g(\CB_{m-1}(f\circ g)+1)\leq g(\CB_{m-1}(f'\circ g')+1)$. 

(2) Proposition 3.2 shows that the function $\CB_{m}^{m}:\FF\times\NN^{m}\ra\NN$ given by the formula $\CB_{m}^{m}(f,b_{1},...,b_{m})=(b_{1}+1)\cdot...\cdot(b_{m}+1)$ has the $m$-bounding property for $m$. Then (1) yields a construction of the function $\CB_{m}^{k}:\FF\times\NN^{k}\ra\NN$ having the $k$-bounding property for $m$ given $\CB_{m-1}:\FF\ra\NN$ (with the bounding property for $m-1$) and $\CB_{m}^{k+1}:\FF\times\NN^{k+1}\ra\NN$ (with the $(k+1)$-bounding property for $m$), for any $k\in\{0,...,m-1\}$. This shows that the function $\CB_{m}^{0}:\FF\ra\NN$ has the bounding property for $m$. \epv

Recall that Proposition 3.1 is the first step of the induction with respect to $m$. The second step of this induction is given in Corollary 3.4 (2). Hence we get the following main result of the section.

\begin{thm} The function $\CB:\NN_{1}\times\FF\ra\NN$ defined recursively in the following way:
\begin{enumerate}[\rm(1)]
	\item $\CB(1,f)=\CB_{1}(f)$ for any $f\in\FF$,
	\item $\CB(m,f)=\CB_{m}^{0}(f)$ for any $m\geq 2$ and $f\in\FF$
\end{enumerate}
has the bounding property.
\end{thm}

{\bf Proof.} It follows from Proposition 3.1 that the function $\CB_{1}:\FF\ra\NN$ such that $\CB_{1}(f)=f(1)+1$ has the bounding property for $m=1$. It follows from Corollary 3.4 (2) that the function $\CB_{m}^{0}:\FF\ra\NN$ has the bounding property for $m$, for any $m\geq 2$. This shows that the construction given in the thesis is correct. \epv

\section{The main results}

In this section we prove the main results of the paper. Throughout we assume that our admissible ordering $\precc$ is graded, e.g. $\precc$ is the degree lexicographical order.

Assume $m\geq 1$ and $f:\N_{1}\ra\N_{1}$ is a function (we do not assume here that $f\in\FF$). An ascending chain $I_{1}\subsetneq ...\subsetneq I_{t}$ of ideals in $K[x_{1},,...,x_{m}]$ is \textit{$f$-bounded} if and only if $I_{j}$ is generated by polynomials of degrees less or equal to $f(j)$, for any $j=1,...,t$. 

Our first goal is to give a bound on the length of $f$-bounded ascending chains of ideals in $K[x_{1},,...,x_{m}]$ depending on $m$ and $f$. The following proposition shows that this problem reduces to the situation studied in Section 3.

\begin{prop} Assume that $m\geq 1$, $f:\N_{1}\ra\N_{1}$ is a function and $I_{1}\subsetneq ...\subsetneq I_{t}$ is an $f$-bounded ascending chain of ideals in $K[x_{1},,...,x_{m}]$. Then there exist monomials $\und{x}^{\alpha_{1}},...,\und{x}^{\alpha_{t}}\in\T_{m}$ such that $\deg(\und{x}^{\alpha_{i}})\leq f(i)$ for $i=1,...,t$ and $\und{x}^{\alpha_{i+1}}\notin\lan \und{x}^{\alpha_{1}},...,\und{x}^{\alpha_{i}}\ran$ for $i=1,...,t-1$. If $f:\N_{1}\ra\N_{1}$ is non-decreasing, this condition is equivalent to the fact that the sequence $\alpha_{1},...,\alpha_{t}$ is an $f$-bounded antichain.
\end{prop}

{\bf{Proof.}} Assume that $I_{j}=\lan h_{j1},h_{j2},...,h_{js_{j}}\ran$ and $\deg(h_{ji})\leq f(j)$ for any $j=1,...,t$ and $i=1,...,s_{j}$. It is easy to see that there are polynomials $h_{1},...,h_{t}$ such that $h_{j}\in\{h_{j1},h_{j2},...,h_{js_{j}}\}$ and $h_{j}\notin\lan h_{1},h_{2},...,h_{j-1}\ran$. Thus we get an ascending chain of ideals of the form $$\lan h_{1}\ran\subsetneq\lan h_{1},h_{2}\ran\subsetneq ...\subsetneq\lan h_{1},h_{2},...,h_{t}\ran$$ of length $t$ with the property that $h_{j+1}\notin\lan h_{1},h_{2},...,h_{j}\ran$ and $\deg(h_{j})\leq f(j)$ for any $j$. We set $H_{1}=\{h_{1}\}$, $H_{2}=\{h_{1},(h_{2})_{H_{1}}\}$, $H_{3}=\{h_{1},(h_{2})_{H_{1}},(h_{3})_{H_{2}}\}$ and so on. We show that $\lan H_{t}\ran=\lan h_{1},h_{2},...,h_{t}\ran$ and the sequence $$\lm(h_{1}),\lm((h_{2})_{H_{1}}),...,\lm((h_{t})_{H_{t-1}})$$ of monomials satisfies the required condition. Indeed, since the admissible ordering $\precc$ is graded, we get $\deg(\lm(h_{1}))=\deg(h_{1})\leq f(1)$ and thus the assertion holds for $t=1$. Assume that the assertion holds for some $t\geq 1$ and there is an ascending chain of ideals  $$\lan h_{1}\ran\subsetneq\lan h_{1},h_{2}\ran\subsetneq ...\subsetneq\lan h_{1},h_{2},...,h_{t}\ran\subsetneq\lan h_{1},h_{2},...,h_{t},h_{t+1}\ran$$ of length $t+1$ such that $h_{j+1}\notin\lan h_{1},h_{2},...,h_{j}\ran$ and $\deg(h_{j})\leq f(j)$ for any $j$. There are polynomials $a_{1},...,a_{t}$ such that $$h_{t+1}=a_{1}h_{1}+a_{2}(h_{2})_{H_{1}}+...+a_{t}(h_{t})_{H_{t-1}}+(h_{t+1})_{H_{t}}$$ and $\lm(h_{t+1})=\max\{\lm(a_{1}h_{1}),\lm(a_{2}(h_{2})_{H_{1}})...,\lm(a_{t}(h_{t})_{H_{t-1}}),\lm((h_{t+1})_{H_{t}})\}$. This yields $h_{t+1}\in\lan H_{t+1}\ran$ and since $\lan h_{1},h_{2},...,h_{t}\ran=\lan H_{t}\ran\subseteq\lan H_{t+1}\ran$, we get that $\lan h_{1},h_{2},...,h_{t},h_{t+1}\ran\subseteq H_{t+1}$. Moreover, $h_{1},(h_{2})_{H_{1}},...,(h_{t+1})_{H_{t}}\in\lan h_{1},h_{2},...,h_{t},h_{t+1}\ran$ and so $\lan H_{t+1}\ran=\lan h_{1},h_{2},...,h_{t},h_{t+1}\ran$. Observe that $(h_{t+1})_{H_{t}}\neq 0$, because otherwise $h_{t+1}\in\lan H_{t}\ran=\lan h_{1},...,h_{t}\ran$, a contradiction. Since $\lm((h_{t+1})_{H_{t}})\precc\lm(h_{t+1})$ and the ordering $\precc$ is graded, we get $\deg(\lm(h_{t+1})_{H_{t}})\leq\deg(\lm(h_{t+1}))\leq f(t+1)$. Finally, the elements of $\LM(H_{t})$ do not divide $\lm((h_{t+1})_{H_{t}})$, because $(h_{t+1})_{H_{t}}$ is reduced modulo $H_{t}$. This implies $\lm((h_{t+1})_{H_{t}})\notin\lan\LM(H_{t})\ran$ which finishes the induction.

To prove the second assertion, assume that $\alpha_{i}=(a_{i1},a_{i2},...,a_{im})\in\NN^{m}$ for $i=1,...,t$. Then $\alpha_{1},...,\alpha_{t}$ is an antichain if and only if for any $i<j$ there is $k$ such that $a_{ik}>a_{jk}$. This implies that the sequence $\und{x}^{\alpha_{1}},...,\und{x}^{\alpha_{t}}$ of monomials in $K[x_{1},...,x_{m}]$ satisfies the conditions $\und{x}^{\alpha_{i+1}}\notin\lan \und{x}^{\alpha_{1}},...,\und{x}^{\alpha_{i}}\ran$ (for $i=1,...,t-1$) and $\deg(\und{x}^{\alpha_{i}})\leq f(i)$ (for $i=1,...,t$) if and only if the sequence $\alpha_{1},...,\alpha_{t}$ is an $f$-bounded antichain. \epv

The above proposition shows that one can associate an $f$-bounded antichain of length $t$ to an $f$-bounded ascending chain of ideals of the same length $t$ (if $f$ is non-decreasing). Therefore we get the following theorem on the length of ascending chains of ideals as a direct consequence of Theorem 3.5 and Proposition 4.1.

\begin{thm} Assume that $m\geq 1$ and $f:\NN_{1}\ra\NN_{1}$ is a function. Suppose that $I_{1}\subsetneq ...\subsetneq I_{t}$ is an $f$-bounded ascending chain of ideals in $K[x_{1},,...,x_{m}]$ of length $t$. Let $g:\NN_{1}\ra\NN_{1}$ be a non-decreasing function such that $g(n)$ is the greatest number of the set $\{f(1),f(2),...,f(n)\}$, for any $n\in\NN$. Then $t\leq\CB(m,g)$. In particular, we have $t\leq\CB(m,f)$, if $f$ is non-decreasing. 
\end{thm}

{\bf Proof.} The chain $I_{1}\subsetneq ...\subsetneq I_{t}$ is $g$-bounded, so the assertion follows from Theorem 3.5 and Proposition 4.1. Note that if $f$ is non-decreasing, then $f=g$. \epv

Now we deduce some consequences of Theorem 4.2 (and hence of Theorem 3.5) in the context of Gr\"obner bases. We start with the following preparatory fact. 



\begin{prop} Assume that $F=\{f_{1},...,f_{s}\}\subseteq K[x_{1},...,x_{m}]$ is a set of non-zero polynomials and $d\geq 1$ is a natural number such that $\deg(f_{i})\leq d$ for $i=1,...,s$. Let $$\lan\LT(B_{0})\ran\subsetneq\lan\LT(B_{1})\ran\subsetneq...$$ be the associated ascending chain of monomial ideals arising from the Buchberger's algorithm. Assume that $n\geq 0$ and $b\in B_{n}$. 
\begin{enumerate}[\rm(1)]
	\item There exist polynomials $a_{1},...,a_{s}\in K[x_{1},...,x_{m}]$ such that $b=a_{1}f_{1}+...+a_{s}f_{s}$ and $\deg(a_{1}),...,\deg(a_{s})\leq (3^{n}-1)d$.
	\item We have $\deg(\lt(b))\leq 3^{n}d$.
\end{enumerate}
\end{prop}

{\bf Proof.} Set $\chi(n)=(3^{n}-1)d$ for any $n\in\NN$. Observe that $(1)$ implies $(2)$. Indeed, if $b=a_{1}f_{1}+...+a_{s}f_{s}$ for some $a_{1},...,a_{s}$ such that $\deg(a_{1}),...,\deg(a_{s})\leq\chi(n)$, then $\deg(a_{1}f_{1}),...,\deg(a_{s}f_{s})\leq\chi(n)+d=3^{n}d$. This implies $\deg(\lt(b))\leq 3^{n}d$, because the ordering $\precc$ is graded.


Thus it is enough to show $(1)$. We use induction with respect to $n$. In the case $n=0$, we have $B_{0}=F$ and $\chi(0)=0$, so the assertion holds. Assume that the assertion holds for some $n\geq 0$, that is, set $B_{n}=\{b_{1},...,b_{r}\}$ and $b_{i}=a_{i1}f_{1}+...+a_{is}f_{s}$ for some $a_{i1},...,a_{is}\in K[x_{1},...,x_{m}]$ such that $\deg(a_{i1}),...,\deg(a_{is})\leq\chi(n)$, for any $i=1,...,r$. We show that the assertion holds for $n+1$. 

Assume that $b_{i},b_{j}\in B_{n}$ and $b_{i}\neq b_{j}$. Recall that $B_{n+1}=B_{n}\cup S_{B_{n}}$ and thus it is enough to show the assertion for $S(b_{i},b_{j})_{B_{n}}$. Observe that $$S(b_{i},b_{j})=\frac{\und{x}^{\alpha}}{\lt(b_{i})}b_{i}-\frac{\und{x}^{\alpha}}{\lt(b_{j})}b_{j}=$$ $$=(\frac{\und{x}^{\alpha}}{\lt(b_{i})}a_{i1}-\frac{\und{x}^{\alpha}}{\lt(b_{j})}a_{j1})f_{1}+...+(\frac{\und{x}^{\alpha}}{\lt(b_{i})}a_{is}-\frac{\und{x}^{\alpha}}{\lt(b_{j})}a_{js})f_{s}$$
where $\und{x}^{\alpha}$ denotes the greatest common multiple of $\lm(b_{i})$ and $\lm(b_{j})$. Since $(1)$ implies $(2)$, we get $$\deg(\frac{\und{x}^{\alpha}}{\lt(b)})\leq\deg(\lt(b'))\leq\chi(n)+d$$ where $b=b_{i}$, $b'=b_{j}$ or vice versa. This yields $$(*)\quad\deg(\frac{\und{x}^{\alpha}}{\lt(b_{i})}a_{ik}-\frac{\und{x}^{\alpha}}{\lt(b_{j})}a_{jk})\leq 2\chi(n)+d,$$ for any $k=1,...,s$, and consequently $\deg(S(b_{i},b_{j}))\leq 2\chi(n)+2d$. Moreover, there are polynomials $c_{1},...,c_{r}$ such that $$S(b_{i},b_{j})_{B_{n}}=S(b_{i},b_{j})-c_{1}b_{1}-...-c_{r}b_{r}=$$ $$=S(b_{i},b_{j})-c_{1}(a_{11}f_{1}+...+a_{1s}f_{s})-...-c_{r}(a_{r1}f_{1}+...+a_{rs}f_{s})$$ and $\lm(c_{t}b_{t})\precc\lm(S(b_{i},b_{j}))$ for any $t=1,...,r$. Because $\precc$ is graded, we get $$\deg(c_{t})\leq\deg(c_{t}b_{t})\leq\deg(S(b_{i},b_{j}))\leq 2\chi(n)+2d$$ and thus $(**)$ $\deg(c_{t}a_{tk})\leq 3\chi(n)+2d$ for any $t=1,...,r$ and $k=1,...,s$. 

It follows by $(*)$ and $(**)$ that the polynomial $S(b_{i},b_{j})_{B_{n}}$ can be written in the form $a_{1}'f_{1}+...+a_{s}'f_{s}$ where $\deg(a_{i}')\leq 3\chi(n)+2d$. Since $3\chi(n)+2d=\chi(n+1)$, this shows the assertion for $n+1$. \epv

By a string $3^{n}d$ we mean the function $f:\NN_{1}\ra\NN_{1}$ such that $f(n)=3^{n}d$ ($d\geq 1$ is a fixed natural number).

\begin{cor} Assume that $F=\{f_{1},...,f_{s}\}\subseteq K[x_{1},...,x_{m}]$ is a set of non-zero polynomials and $d$ is a natural number such that $\deg(f_{i})\leq d$ for $i=1,...,s$. Let $$\lan\LT(B_{0})\ran\subsetneq\lan\LT(B_{1})\ran\subsetneq...\subsetneq\lan\LT(B_{r})\ran$$ be the associated ascending chain of monomial ideals arising from the Buchberger's algorithm such that $B_{r}$ is the Gr\"obner basis of $\lan F\ran$. Then $r+1\leq\CB(m,3^{n}d)$.
\end{cor}

{\bf Proof.} It follows from Proposition 4.3 (2) that the ascending chain $\lan\LT(B_{0})\ran\subsetneq\lan\LT(B_{1})\ran\subsetneq...\subsetneq\lan\LT(B_{r})\ran$ is $3^{n}d$-bounded. Hence Theorem 4.2 yields the condition $r+1\leq\CB(m,3^{n}d)$. \epv

Set $m\geq 1$, $d\geq 1$ and define the function $\gamma_{m,d}:\NN\ra\NN$ in the following way $$\gamma_{m,d}(i)=(3^{\CB(m,3^{n}d)-1}-1)d+i$$ for any $i\in\NN$. The function $\gamma_{m,d}$ has the following property.

\begin{cor} Assume that $m\geq 1$ and $d\geq 1$. Then for any $g\in K[x_{1},...,x_{m}]$ and $f_{1},...,f_{s}\in K[x_{1},...,x_{m}]$ such that $\deg(f_{i})\leq d$ for $i=1,...,s$ the following condition is satisfied: $g\in\lan f_{1},...,f_{s}\ran$ if and only if there exist $h_{1},...,h_{s}\in K[x_{1},...,x_{m}]$ such that $g=h_{1}f_{1}+...+h_{s}f_{s}$ and $\deg(h_{i})\leq\gamma_{m,d}(\deg(g))$ for $i=1,...,s$.
\end{cor}

{\bf Proof.} Assume that $g,f_{1},...,f_{s}\in K[x_{1},...,x_{m}]$ and $\deg(f_{i})\leq d$ for $i=1,...,s$. Set $F=\{f_{1},...,f_{s}\}$ and let $\lan\LT(B_{0})\ran\subsetneq\lan\LT(B_{1})\ran\subsetneq...\subsetneq\lan\LT(B_{r})\ran$ be the ascending chain of monomial ideals arising from the Buchberger's algorithm such that $B_{r}=G=\{g_{1},...,g_{t}\}$ is the Gr\"obner basis of $\lan F\ran$. 

Assume that $g\in\lan f_{1},...,f_{s}\ran$. Since $G$ is a Gr\"obner basis of $\lan F\ran$, there are polynomials $p_{1},...,p_{t}$ such that $g=p_{1}g_{1}+...+p_{t}g_{t}$ and $\lm(g)$ is the maximal element of $\{\lm(p_{1}g_{1}),...,\lm(p_{t}g_{t})\}$. Hence $\lm(p_{i}g_{i})\precc\lm(g)$ so $\deg(p_{i}g_{i})\leq\deg(g)$ and consequently $\deg(p_{i})\leq\deg(g)$, for any $i=1,...,t$.

Corollary 4.4 yields $r+1\leq\CB(m,3^{n}d)$. Furthermore, Proposition 4.3 (1) implies that $g_{i}=a_{i1}f_{1}+...+a_{is}f_{s}$ for some polynomials $a_{i1},...,a_{is}$ with $$\deg(a_{i1}),...,\deg(a_{is})\leq(3^{r}-1)d\leq(3^{\CB(m,3^{n}d)-1}-1)d,$$ for $i=1,...,t$. It follows that $$\deg(p_{i}a_{ik})\leq (3^{\CB(m,3^{n}d)-1}-1)d+\deg(g)=\gamma_{m,d}(\deg(g))$$ for $i=1,...,t$ and $k=1,...,s$. This shows the assertion. \epv

Let us note that the main results of this section (Theorem 4.2 and Corollaries 4.4 and 4.5) do not depend on the choice of the base field $K$ of the polynomial ring $K[x_{1},...,x_{m}]$. 

\section{Remarks}

Our motivation to study problems concerning ascending chains of ideals in the polynomial ring arises from the first order logic. The goal is to give a constructive proof of the renowned Tarski's theorem on quantifier elimination in the theory of algebraically closed fields. This theorem was proved by A. Tarski in 1948 in an unpublished paper, see \cite{R} for the details.

Roughly, Tarski's theorem states that if $\varphi(x_{1},...,x_{n})$ is a formula in the first order language of the theory of fields with $n$ free variables $x_{1},...,x_{n}$, then there exists a quantifier-free formula $\varphi'(x_{1},...,x_{n})$ (a formula in which quantifiers do not occur), with the same free variables, such that $\varphi(x_{1},...,x_{n})$ is equivalent with $\varphi'(x_{1},...,x_{n})$. This means that for any algebraically closed field $K$ and any elements $a_{1},...,a_{n}\in K$ we have $\varphi(a_{1},...,a_{n})\leftrightarrow\varphi'(a_{1},...,a_{n})$. We refer to \cite{M} for the necessary details. 

As an example, consider the formula $\varphi(A)=\exists_{B\textnormal{ }}AB=BA=I_{n}$ where $A,B$ are $n\times n$ complex matrices and $I_{n}$ is the $n\times n$ identity matrix ($\varphi(A)$ can be suitably written in the first order language of the theory of fields). This formula states that $A$ is non-singular and thus $\varphi(A)$ holds if and only if $\det(A)\neq 0$. The latter formula is quantifier-free and very easy to verify. Generally, this is the case for any quantifier-free formula.




Standard proofs of Tarski's theorem are existential, that is, they do not provide the form of the quantifier-free formula equivalent with the given one. A constructive proof aims to provide that form. In the subsequent paper \cite{P} we apply Corollary 4.5 to give a constructive proof of Tarski's theorem. Moreover, we show some interesting applications of this constructive version. For example, a formula stating the existence of a common invariant subspace of $n\times n$ complex matrices $A_{1},...,A_{s}$ is a first order formula $\psi$ of the theory of fields. By the constructive Tarski's theorem we are able to give a quantifier-free formula $\psi'$ which is equivalent to $\psi$. The formula $\psi'$ may be considered as an algorithm for verifing the existence of a common invariant subspace of $A_{1},...,A_{s}$. We emphasize that until \cite{ArPe}, published in 2004, it was not known if such an algorithm exists in the general case (for special cases see \cite{Sh}, \cite{Ts} and \cite{JP}). 

In the series of papers \cite{JKP}, \cite{JP}, \cite{PJ} and \cite{PKJ} we consider algorithms (called there \textit{computable conditions}) for the existence of various common invariant subspaces of complex linear operators. We further apply these algorithms in some problems of quantum information theory. All the problems we consider can be expressed in the first order language of the theory of fields, and hence the constructive Tarski's theorem is applicable. This gives a new general context for this research and opens the possibility for other applications. Note that some impact of quantifier elimination technique on quantum information theory has been recently noticed in \cite{WCP-G}.

\section*{Acknowledgements} This research has been supported by grant No. DEC-2011/02/A/ST1/00208 of National Science Center of Poland.

\noindent Grzegorz Pastuszak\\Faculty of Mathematics and Computer Science\\ Nicholaus
Copernicus University\\ Chopina 12/18\\ 87-100 Toru\'n, Poland\\
past@mat.uni.torun.pl

\end{document}